\documentclass{article}%
\usepackage{amsmath}
\usepackage{amsfonts}
\usepackage{amssymb}
\usepackage{graphicx}%
\newtheorem{theorem}{Theorem}

\newtheorem{lemma}[theorem]{Lemma}

\begin{document}

\title{Note on the classification of the orientation reversing homeomorphisms of
finite order of surfaces}
\author{Antonio F. Costa \thanks{Partially supported by PGC2018-096454-B-100 (Spanish
Ministry of Science, Innovation and Universities.)}}
\date{}
\maketitle

\section{Introduction}

The classification of periodic orientation reversing autohomeomorphims of a
closed oriented surface has been made by Kazuo Yokoyama, Complete
classification of periodic maps on compact surfaces, Tokyo J. Math. 15 (1992)
and A. F. Costa, Classification of the orientation reversing homeomorphisms of
finite order of surfaces, Topology and its Applications 62 (1995) following
different approaches. Weibiao Wang of Peking University and Chao Wang of East
China Normal University in Shanghai, have pointed out certain errors in
Theorems 4.4 and 4.5 of the last article (in case of orientation reversing
homeomorphims of period multiple of 4). These errors make it is necessary to
modify the definition of the invariant $h_{1}$ and replace the condition (iii)
of the main Theorem 0.2 with condition (ii) of Theorem 1 as shown in Section
3. In the Note of Section 3 we present the necessary changes in the statements
of the Theorems 4.4 and 4.5 and a correction on the definition of isotropy
invariants for the case of orientable quotient surfaces. In Section 4 we
correct in the wording of propositions of some subsequent works.

I thank very much to Weibiao Wang and Chao Wang for the messages and
communications on this subject and to the referee for careful reading and
several suggestions and corrections.

\bigskip

\bigskip

\section{Preliminaries and notations}

Let $S$ be an orientable closed surface and $\Phi$ be an orientation reversing
autohomeomorphism of finite order $2q$ of $S$. Two homeomorphisms $\Phi_{1}$
and $\Phi_{2}$ of a surface $S$ are said to be topologically equivalent if
there is a homeomorphism $h:S\rightarrow S$ such that $\Phi_{1}=h^{-1}%
\circ\Phi_{2}\circ h$.

The orbifold structure on the orbit space will be denoted by $S/\Phi$ and
$\left\vert S/\Phi\right\vert $ is the underlying quotient surface. The
projection $\pi:S\rightarrow S/\Phi$ is an orbifold covering.

We assume $q$ is even and greater than $2$, then $\left\vert S/\Phi\right\vert
$ is a closed non-orientable surface. The orbifold covering $\pi:S\rightarrow
S/\Phi$ has a finite set of singular values corresponding to the conic points
of the orbifold $S/\Phi$. Let $r$ be the number of conic points of $S/\Phi$.
There are canonical presentations of the orbifold fundamental group $\pi
_{1}O(S/\Phi)$ as follows:%
\[
\left\langle d_{1},...,d_{g},x_{1},...,x_{r}:x_{1}...x_{r}d_{1}^{2}%
...d_{g}^{2}=1,x_{i}^{m_{i}}=1,i=1,...,r\right\rangle
\]
where $g$ is the topological genus of $\left\vert S/\Phi\right\vert $. The
relation $x_{1}...x_{r}d_{1}^{2}...d_{g}^{2}=1$ is called the long relation.
The canonical generators $x_{i}$ of two canonical presentations are conjugate
or inversed.

For the abelianization $H_{1}O(S/\Phi)$ of $\pi_{1}O(S/\Phi)$, we define
\emph{canonical generator system} to be the set of generators obtained from a
canonical presentation of $\pi_{1}O(S/\Phi)$. The canonical generators of
$H_{1}O(S/\Phi)$ will be denoted by
\[
X_{1},...,X_{r},D_{1},...,D_{g}.
\]
The capital letter denotes the homology class determined in $H_{1}O(S/\Phi)$
by the generator of $\pi_{1}O(S/\Phi)$ denoted by the corresponding small
letter. From the long relation we have the relation $2D_{1}+...+2D_{g}%
+\sum_{i=1}^{r}X_{i}=0$.

We choose a non-singular base point $o\in S/\Phi$ and $\widetilde{o}\in
\pi^{-1}(o)$ and we name by $o_{i}=\Phi^{i}(\widetilde{o})$, the elements of
$\pi^{-1}(o)$. We take as generator of $Z_{2q}=\left\langle [1]\right\rangle $
the cycle $(1,2,...,2q)=[1]$. The covering $\pi:S\rightarrow S/\Phi$
determines an epimorphism $\pi_{1}O(S/\Phi,o)\rightarrow Z_{2q}\cong%
\left\langle (1,2,...,2q)\right\rangle \subset\mathrm{Sym}\{1,...,2q\}$, and
since $Z_{2q}$ is abelian we have the monodromy epimorphism $T:H_{1}%
O(S/\Phi)\rightarrow Z_{2q}$.

If $T_{1}$ is the monodromy of $\Phi_{1}$ and $T_{2}$ is the monodromy of
$\Phi_{2}$ the homeomorphisms $\Phi_{1}$ and $\Phi_{2}$ are topologically
equivalent if there is an automorphism $\alpha$ of $\pi_{1}O(S/\Phi)$ inducing
an automorphism $\alpha_{\ast}$ of $H_{1}O(S/\Phi)$ such that $T_{1}%
=T_{2}\circ\alpha_{\ast}$.

Since each automorphism $\alpha_{\ast}$ sends canonical presentations to
canonical presentations we have that the set $\{\pm T(X_{i}),$ $i=1,...,r\}$
is a topological invariant that we shall call the \emph{set of isotropies}.

We shall use the following automorphisms $H_{1},...,H_{4}$ of $H_{1}O(S/\Phi)$
(induced by automorphisms of $\pi_{1}O(S/\Phi)$):

- $H_{1}(i,j)(D_{i})=2D_{j}+D_{i};\,H_{1}(i,j)(D_{j})=-D_{j},$ $(1\leq i\leq
g,1\leq j\leq g,i\neq j)$

- $H_{2}(i,j)(D_{i})=D_{i}+X_{j};\,H_{2}(i,j)(X_{j})=-X_{j},$ $(1\leq i\leq
g,1\leq j\leq r)$

- $H_{3}(i,j)(X_{i})=X_{j};\,H_{3}(i,j)(X_{j})=X_{i}$, if order$(X_{i})=$
order$(X_{j}),$ $(1\leq i\leq r,1\leq j\leq r)$

- $H_{4}(i,j)(D_{i})=D_{j};\,H_{4}(i,j)(D_{j})=D_{i}$, $(1\leq i\leq g,1\leq
j\leq g).$

The remaining generators of $H_{1}O(S/\Phi)$ that do not appear in the above
formulae are unchanged by the automorphisms. $H_{1},H_{2},H_{3}$ correspond,
respectively, to $H_{1}^{\ast},H_{2}^{\ast},H_{3}^{\ast}$ [1, pages 152-154],
and $H_{4}$ arises from a pants move [1, page 149], when the two swapped
boundary components bound Mobius strips.

\section{Classification of orientation reversing autohomeomorphims of period a
multiple of 4}

In order to establish the classification we need to define two topological
invariants $h_{1}$ and $h_{2}$, the definition is very close to the given in
\cite{Top} but we need to make some essential changes.

We have an orientable closed surface $S$ and an orientation reversing
autohomeomorphism $\Phi$ of finite order $2q$ of $S$.

The invariant $h_{1}$ is defined for orientation reversing homeomorphisms
without isotropies of type $T(X_{i})=[q]$.

Let $D_{1},...,D_{g},X_{1},...,X_{r}$ be a canonical generator system of
$H_{1}O(S/\Phi)$. Using the automorphism $H_{2}$ of section 2 we may modify
the generators such that satisfy: $T(X_{i})\in\{[2],...,[q-2]\}$. Using the
relation $2D_{1}+...+2D_{g}+\sum_{i=1}^{r}X_{i}=0$, the element $T(D_{1}%
+...+D_{g})$ is determined by the $T(X_{i})$ up addition by $[q]$. We define
the invariant $h_{1}(\Phi)=0$ if $T(D_{1}+...+D_{g})\in\{[0],...,[q-1]\}$ and
$h_{1}(\Phi)=1$ if $T(D_{1}+...+D_{g})\in\{[q],...,[2q-1]\}$.

The geometrical interpretation of this invariant is given by the monodromy of
a homology class $c$ of $H_{1}O(S/\Phi)$ that is represented by a closed curve
$\gamma$ such that cutting $\left\vert S/\Phi\right\vert $ through $\gamma$ we
obtain an orientable surface and the position of this curve with respect to
the conical points of $S/\Phi$ is given by the condition $T(X_{i}%
)\in\{[2],...,[q-2]\}$.

The invariant $h_{2}$ is necessary only in the case where $\left\vert
S/\Phi\right\vert $ has genus two. Let $l$ be the smallest integer such that
the element $[2l]$ is a generator of the subgroup of $Z_{2q}$ generated by the
isotropies. Let $\left\vert S/\Phi^{2l}\right\vert $ be the orbit surface of
the action of $\Phi^{2l}$. The homeomorphism $\Phi$ defines an orientation
reversing finite order and fixed point free homeomorphism $\Phi^{free}$ on
$\left\vert S/\Phi^{2l}\right\vert $. The invariant $h_{2}$ is given by the
topological type of $\Phi^{free}$. The first homology $H_{1}(\left\vert
S/\Phi\right\vert )$ is isomorphic to $Z\oplus Z_{2}=\left\langle
x,y:x+y=y+x,2y=0\right\rangle $ and all the automorphisms of $Z\oplus Z_{2}$
send $x$ to an element in the set $S=\{\delta x+\varepsilon y:\delta
=-1,1;\varepsilon=0,1\}$. If we note by $T_{free}:H_{1}(\left\vert
S/\Phi\right\vert )\rightarrow Z_{2q}$ the monodromy of $\Phi^{free}$, the set
$T_{free}(S)$ is a topological invariant that we denote it by $h_{2}$:
\begin{gather*}
h_{2}=\{\delta T_{free}(D_{1})+\varepsilon T_{free}(D_{1}+D_{2}):\delta
=-1,1;\varepsilon=1,0\}=\\
\{\delta T_{free}(D_{1})+\varepsilon h_{1}(\Phi^{free}):\delta
=-1,1;\varepsilon=1,0\}.
\end{gather*}

With this invariants, it is now possible to establish a classification theorem
that corrects the Theorem 0.2 in the introduction of \cite{Top}:

\begin{theorem}
Let $\Phi_{1}$ and $\Phi_{2}$ be two orientation reversing autohomeomorphisms
of finite order $2q$ of a surface $S$. Assume $q$ is even.

The homeomorphisms $\Phi_{1}$ and $\Phi_{2}$ are topologically equivalent if
and only if the following three statements are all true:

(i) $\Phi_{1}$ and $\Phi_{2}$ have the same set of isotropies.

(ii) if there is not any isotropy of order $2$, $h_{1}(\Phi_{1})=h_{1}%
(\Phi_{2})$.

(iii) if $\left\vert S/\Phi\right\vert $ has genus two and $\Phi_{1}^{free}$,
$\Phi_{2}^{free}$ have order greater than two (equivalently the set of
isotropies is not a generator system of $\left\langle [2]\right\rangle \leq
Z_{2q}$), $h_{2}(\Phi_{1})=h_{2}(\Phi_{2})$.
\end{theorem}

\bigskip

\textbf{Proof.}

The set of isotropies, $h_{1}$ and $h_{2}$ are topological invariants by the
very definitions, so it suffices to prove that these invariants determine the
topological type. The way of proving that is checking that these invariants
determine completely the monodromy $T$ of a given orientation reversing
autohomeomorphism $\Phi$.

We consider $S/\Phi$ of genus $g$ and with $r$ conical points. Let
\[
D_{1},...,D_{g},X_{1},...,X_{r}%
\]
be a canonical system of generators of $H_{1}O(S/\Phi)$.

\textbf{Case 1}. Genus $g$ of $\left\vert S/\Phi\right\vert $ different from two.

We shall use the following Lemma (see Lemma 3.1(2) of \cite{Top}):

\begin{lemma}
If $\left\vert S/\Phi\right\vert $ has genus $g>2$ then there is an
automorphism $h$ of $S/\Phi$ such that: $T(h(D_{1}))=...=T(h(D_{g-1})=[1]$,
$h(X_{i})=X_{i}$, $i=1,...,r$.
\end{lemma}

Using the Lemma we can assume $T(D_{i})=[1],i\neq g$.

Subcase 1. There is at least one isotropy of order 2.

By automorphisms $H_{3}(i,j)$ we can assume that $T(X_{1})=...=T(X_{s}%
)=[q],$for some $s\geq1$ and using automorphisms $H_{2}(i,j)$ we can have
$T(X_{i})\in\{[2],...,[q-2]\}$ $i>s$.

By the long relation $2D_{1}+...+2D_{g}+\sum_{i=1}^{r}X_{i}=0$ and
$H_{2}(i,1)$, we can have $T(D_{g})\in\{[1],...,[q-1]\}$ and this fact
determines $T(D_{g})$ and then completely $T$.

Subcase 2. There is no $T(X_{j})=[q]$.

By automorphisms $H_{2}(i,j)$ we can assume $T(X_{i})\in\{[2],...,[q-2]\}$.
Then using the long relation $2D_{1}+...+2D_{g}+\sum_{i=1}^{r}X_{i}=0$ and
invariant $h_{1}$ we determine completely $T(D_{g})$ and the topological type
of $\Phi$.

\bigskip

\textbf{Case 2}. $\left\vert S/\Phi\right\vert $ has genus two.

By automorphisms $H_{2}$ and $H_{4}$ we can assume $T(X_{i})=[q]$ and
$T(X_{j})\in\{[2],...,[q-2]\}$, $i=1,...,s$, $j=s+1,...,r$, and using the
isotropies of order $2$ or the invariant $h_{1}$, as in the previous case, we
can determine $T(D_{1}+D_{2})$. Let $l$ be the smallest integer such that the
element $[2l]$ is a generator of the subgroup of $Z_{2q}$ generated by the
isotropy invariants. By the automorphisms $H_{2}(1,j)$ and $H_{2}(2,j)$ we can
obtain $T(D_{1})\in\{[1],...,[l]\}$ and using if necessary $H_{1}(1,2)$ and
$h_{2}$ the value of $T(D_{1})$ is determined. \ $\square$

\bigskip

\textbf{Note.} In \cite{Top}, Theorems 4.4 and 4.5 are wrong as stated. The
hypothesis in Theorem 4.4: \textquotedblleft$(\Phi_{j}/4),$ $j=1,2$, have
fixed points\textquotedblright\ must be replaced by \textquotedblleft there is
some $T_{j}(X_{i}^{(j)})=[q]$\textquotedblright. In Theorems 4.5 and 4.6 the
hypothesis \textquotedblleft$(\Phi_{j}/4)^{2},$ $j=1,2$, are fixed point
free\textquotedblright\ must be replaced by \textquotedblleft there is no any
$T_{j}(X_{i}^{(j)})=[q]$\textquotedblright. In Theorems 4.5 and 4.6 the
invariants $h_{1}$ and $h_{2}$ must be as defined above in this Section.

Another correction to \cite{Top} is that, in the case (not previously
considered here) where $S/\Phi$ is orientable with $s\geq0$ boundary
components, the invariant given by the set of isotropies is the unordered pair
of sets
\begin{align*}
\{T(X_{i}),T(C_{j})  &  :i=1,...,r;\text{ }j=1,...,s\},\\
\{-T(X_{i}),-T(C_{j})  &  :i=1,...,r;\text{ }j=1,...,s\}
\end{align*}
where $r$ is the number of conic points in $S/\Phi$ and $s$ in the number of
boundary components of $\left\vert S/\Phi\right\vert $. Note that the
orientation of $\left\vert S/\Phi\right\vert $ is not determined by the
orientation of $S$.

\bigskip

\section{Consequences and corrections in \cite{Cont} and \cite{Comm}}

The classification in \cite{Top} has been used in \cite{Cont} and in
\cite{Comm}. We correct here the statements following Section 3.

In the Proposition 2 in \cite{Cont} the hypothesis $f^{q}$ have fixed points
must be replaced by there are branched points with isotropy groups of order 2.

In \cite{Comm} it is necessary to change the statement of Proposition 2.1
which establishes the topological types of the anticonformal automorphisms of
order a multiple of $4$ that can be represented as a symmetry of an embedded
surface in Euclidean space (embeddable automorphism). Let $f$ be an
anticonformal automorphism of order $q$ with $q$ even of a closed Riemann
surface $S$ (note that in this case $S/f$ is non orientable, with empty
boundary and genus $g$). Let $D_{1},...,D_{g},X_{1},...,X_{r}$ be a canonical
generator system of $H_{1}O(S/f)$, and $T:H_{1}O(S/f)\rightarrow Z_{2q}$ be
the monodromy of the covering $S\rightarrow S/f$. The statement of Proposition
1.2 must be as follows:

\begin{theorem}
The anticonformal automorphism $f$ of order $2q$, $q$ even, is embeddable if
and only if it satisfies one of the following conditions:

1. $f^{2}$ is fixed point free and $h_{1}(f)=0$,

2. $f^{2}$ has fixed points and $q=2,$ $g\geq r,$

3. $f^{2}$ is has $r$ fixed points, $r>0$, $q>2$, $g\geq r$, $r+g\equiv
0\operatorname{mod}2$, all the $T(X_{i})$ are equal to $[\pm m]$, where
$[m]\in\{[2],...,[q-2]\}$ is some generator or $\left\langle [2]\right\rangle
$, and if $[\frac{rm}{2}]\in\{[1],...,[q-1]\}$ the invariant $h_{1}(f)$ must
be equal $1$, if $[\frac{rm}{2}]\in\{[q+1],...,[2q-1]\}$ the invariant
$h_{1}(f)$ must be $0$.
\end{theorem}

\end{document}